# A novel reduced basis method for adjoint sensitivity analysis of dynamic topology optimization

Shuhao Li[a,b], Hu Wang[a,b*], Jichao Yin[a,b], Xinchao Jiang[a,b], Yaya Zhang[a,b]

[a]*State Key Laboratory of Advanced Design and Manufacturing for Vehicle Body, Hunan University, Changsha 410082, People's Republic of China;*

[b]*Beijing Institute of Technology Shenzhen Automotive Research Institute, Shenzhen 518000, People's Republic of China*

[*]Hu Wang: wanghu@ szari.ac.cn

## Highlights

- An efficient RBM-based adjoint sensitivity analysis method is proposed
- A novel model-based error estimation in reduced basis method is developed
- An offline basis vector generation strategy for reduced-order model is advanced
- An adaptive strategy for the online phase in reduced basis method is introduced

## Abstract

In gradient-based time domain topology optimization, design sensitivity analysis (DSA) of the dynamic response is essential, and requires high computational cost to directly differentiate, especially for high-order dynamic system. To address this issue, this study develops an efficient reduced basis method (RBM)-based discrete adjoint sensitivity analysis method, which on the one hand significantly improves the efficiency of sensitivity analysis and on the other hand avoids the consistency errors caused by the continuum method. In this algorithm, the basis functions of the adjoint problem are constructed in the offline phase based on the greedy-POD method, and a novel model-based estimation is developed to facilitate the acceleration of this process. Based on these basis functions, a fast and reasonably accurate model is then built by Galerkin projection for sensitivity analysis in each dynamic topology optimization iteration. Finally, the effectiveness of the error measures, the efficiency and the accuracy of the presented reduced-order method are verified by 2D and 3D dynamic structure studies.

**Keywords** Reduced basis method; Adjoint sensitivity analysis; Error estimation; Dynamic topology optimization

## 1. Introduction

Gradient-based topology optimization typically involves a large number of quantifiable design variables [1], making an efficient sensitivity analysis method crucial. Sensitivity analysis methods can be divided into three categories regarding structural performance: approximate, discrete and continuum [2]. Approximate sensitivity analysis is mainly applied to finite difference or complex-step derivative [3, 4]. The core of the discrete method is to first discretize the structural governing equation and then perform the differential operation, also known as the discretize-then-differentiate method, whereas the continuum method, referred to as the differentiate-then-discretize method, lies in taking the design derivative of the variational equation for continuum model [2, 5]. We adopt the discrete approach for adjoint sensitivity analysis to eliminate consistency errors observed in the commonly-used continuum approach [6]. Furthermore, from the point of view of the sensitivity solution, there can be a distinction between the direct method and the adjoint method. The direct method involves directly solving the partial derivative of the state variable with respect to the design variable, and then obtaining the sensitivity solution by differentiating the chain expression. As for the adjoint method, the original problem is transformed into a computationally tractable adjoint problem by introducing adjoint variables that contain all implicitly dependent terms, thus avoiding the direct derivation of state variables for each design variable. In particular, the adjoint method has advantages in the treatment of structure optimization problem with a large number of design variables [2, 7-10]. For topology optimization problem, the design variables are typically defined at the element level. The number of design variables is much larger than the number of responses, so the adjoint method is more efficient than direct method in terms of sensitivity solution [11]. For the topology optimization of transient problems in this study, the adjoint method requires an inverse transient analysis to obtain the adjoint variable. As the size of the dynamic system increases in space and time, the cost of computing the adjoint sensitivity increases significantly [12].

In order to avoid the so-called curse of dimensionality in the above adjoint sensitivity calculation for dynamic topology optimization, a model order reduction approach should be used. Model order reduction (MOR) aims to transform high-dimensional problems into low-dimensional space using a specific feature extraction technique, while maintaining relatively high approximation quality. At present, these MOR techniques, including reduced basis method [13], Ritz vector method [14], mode displacement method [15], etc., have been widely used in time-domain dynamic response calculation and design dynamic sensitivity analysis due to their high efficiency and low storage requirements [16].

The majority of MOR-based studies mentioned above use the direct method for the dynamic sensitivity analysis. For the adjoint method, based on Xiao's studies on static topology optimization with discrete adjoint [17], Qian[18] further developed the on-the-fly reduction method for dynamic topology optimization problem with continuum adjoint. In addition, according to the author's review of the relevant literature, most studies on the time-domain dynamic topology optimization are based on continuum adjoint approach [19-21], which suffers from consistency error in sensitivity analysis. To address this issue, this study proposes a novel discrete adjoint sensitivity analysis method that integrates an improved reduced basis method to project high-dimensional adjoint problems into low-dimensional spaces, significantly improving computational efficiency while maintaining sensitivity analysis accuracy.

The present reduced basis method (RBM) consists of an offline and an online phase. In the offline phase, the features of the basic model are extracted by acquiring a typically small set of basis functions that are problem specific. The basic model is represented by snapshots obtained by computing the solution of the full-order model (exact solution) for a representative parameter space. The basis functions can be obtained by various methods: the most widely used are proper orthogonal decomposition (POD) and the greedy algorithm [22, 23]. With the basis functions, an affine relevance of the full-order model on the parameter space is exploited to build a fast-to-compute low-order model. In the online phase, the approximate solution can be solved for each new set of parameters by solving the low-order model. In general, the dimension of the low-order model is much smaller than that of the full-order model, requiring less computation.

Besides the construction of the reduced-order model, it is crucial for any MOR technique to assess and determine the accuracy of the reduced model, which can be achieved by error estimation. Thus, a fast and relatively accurate posteriori error estimation is indispensable in RBM. First, it supports the greedy algorithm to speed up the generation of basis vectors in the offline. And the second role is to assess the accuracy of the approximate solution generated by solving the low-order model in the online. Numerous studies have been carried out to develop a residual-based error estimation for elliptic PDEs [24, 25], hyperbolic [26, 27] and parabolic [23, 28]. Haasdonk [29] derived an error bound from the norm of the state matrix, but it increases exponentially with time if the norm is greater than 1, so it cannot accurately evaluate the reduced order model. Zhang [27] developed a novel error estimation supported by a dual system, which effectively avoids the accumulation of residual over time. Bernreuther et al. [24] proposed an efficient error estimation for the adjoint problem that avoids overestimating the true adjoint error. Xu [25] developed to estimate the error for second order elliptic problems in accordance with the Helmholtz decomposition yield, which is still appropriate for certain problems. To this end, Abbasi et al. [26] proposed an efficient and simple error bound in terms of the Euclidean distance, which does not depend on the norm of the state matrix, but this error estimation relies too much on offline sampling to obtain the exact estimate error. Furthermore, the error bound would be approximate owing to its non-mathematical nature. The ultimate problem is that there are rarely efficient error estimators that are both versatile and accurate.

Inspired by Abbasi above, further study is conducted. A model-based error estimation method is proposed to establish the relationship between the residual error and the true error. In recent decades, Natural Network (NN) models have been widely applied in many fields, such as prediction and classification, etc., due to their excellent modelling capabilities [30, 31]. The Fully-Connected Neural Network (FNN), as one of the most representative among NNs, is a multi-layer feedforward neural network, which can realize the complex relationship mapping of input and output and self-learning. Therefore, this study constructs an FNN to estimate the true adjoint error in dynamic optimization.

As a matter of fact, some scholars have attempted to introduce NN techniques such as Convolutional Neural Networks (CNN) [32], Generative Adversarial Network (GAN) [33] etc. into topology optimization to speed up the optimization process. In fact, these studies lie in replacing the optimizer with NN models to efficiently obtain the final designs through a train-then-predict process, where the initial or partially converged designs are mapped to the final designs. Unlike the above applications, in this study, the FNN is used to construct error estimators in a reduced-order model for dynamic sensitivity analysis to accelerate topology optimization process. To the best of the author's knowledge, there is no relevant research on the application of NN to error estimation in reduced-order models and also to dynamic topology optimization. The main contributions of this study are described below.

In this study, we propose an efficient strategy for dealing with sensitivity analysis in dynamic topology optimization. A discrete adjoint sensitivity analysis is adopted to eliminate consistency errors. Furthermore, a reduced basis method associated with this adjoint sensitivity analysis for dynamic topology optimization is proposed and tested on an elliptic system. Finally, a novel model-based error estimator is developed for this proposed reduced-order model.

The outline of this study is as follows. In Section 2, the general dynamic topology optimization with the discrete adjoint sensitivity analysis scheme is explained and the full-order model of the adjoint problem is derived. In Section 3, the ingredients of the reduced basis method and the reduced-order model are introduced. In Section 4, the existing error estimation methods are discussed and the model-based error estimator is developed. In Section 5, the results of the three numerical examples under different external load conditions are presented. Finally, the conclusions and suggestions are given.

## 2. Problem statement

This section introduces the theoretical aspects of dynamic topology optimization and presents the main problems to be studied. In this study, the following governing dynamic equilibrium equations are studied considering Rayleigh damping for a multi-DOF:

$$\mathbf{M}\ddot{\mathbf{d}}(t) + \mathbf{C}\dot{\mathbf{d}}(t) + \mathbf{K}\mathbf{d}(t) = \mathbf{F}(t) \tag{1}$$

where $\mathbf{M}$, $\mathbf{C}$, $\mathbf{K}$ denote the mass, damping and stiffness matrices respectively. The unknown displacement, velocity, acceleration and external force vectors are functions of the time variable $t$. The solution of the unknown response is discussed in the following

subsection. The Rayleigh damping using proportional damping is expressed in a discrete form as follows:

$$\mathbf{C} = \alpha_M \mathbf{M} + \alpha_K \mathbf{K} \qquad (2)$$

in which $\alpha_M$ and $\alpha_K$ are the Rayleigh damping parameters.

In Section 2.1, we present the topology optimization framework for linear elastic structures under dynamic loading by a continuous density field. In Section 2.2, we present the discrete adjoint sensitivity analysis method and derive the reduced sensitivity analysis strategy.

## *2.1. Dynamical topology optimization problem*

Topology optimization aims to find the material distribution of the structure by optimizing a given objective function containing some performance metrics while satisfying the constraint requirements. In the density-based topology optimization method, the element density variable ($b$) is defined to describe the material distribution of the structure, which covers the range of $[0,1]$, from no material to the available material. In general, the objective function $f(b,\mathbf{d})$ and the constraints $g_i(b,\mathbf{d})$, $i = 1,...,l$, depend on the density variables, $b$, and physical fields, $\mathbf{d}$, such as displacement, velocity, stress, etc., which are obtained by solving the governing equations of specific problems. The above density-based topology optimization problem can be written in formal terms as follows:

$$\inf_{b \in \xi} f(b,\mathbf{d})$$
$$\text{s.t.} \quad g_i(b,\mathbf{d}) \leq 0, \qquad i = 1,...,L \qquad (3)$$

where $\xi$ is the space of prescribed density variables, $L$ the number of constraints.

Given that the subject of this study is the elastodynamic problem described in Eq. (1), the involved state variables are defined such that, given an initial displacement, $\mathbf{d}_0(\mathbf{x})$, initial velocity, $\mathbf{v}_0(\mathbf{x})$, given displacements, $\bar{\mathbf{d}}$ (applied to the portion, $\Gamma_D$, of $\partial \Omega$), and given boundary tractions, $\mathbf{F}$ (applied to the portion, $\Gamma_F$, of $\partial \Omega$), the displacement field, $\mathbf{d} \in \zeta$, is obtained by the following boundary value problem:

$$\begin{aligned}
\text{div } \boldsymbol{\sigma} + \mathbf{F}_b = \bar{m}\ddot{\mathbf{d}} + \bar{c}\dot{\mathbf{d}}, & \quad \text{in } \Omega \times (0,T] \\
\mathbf{d}(\mathbf{x},t) = \bar{\mathbf{d}}(\mathbf{x},t), & \quad \text{on } \Gamma_D \times (0,T] \\
\boldsymbol{\sigma}(\mathbf{x},t) \cdot \mathbf{n} = \mathbf{F}(\mathbf{x},t), & \quad \text{on } \Gamma_F \times (0,T] \\
\mathbf{d}(\mathbf{x},0) = \mathbf{d}_0(\mathbf{x}), & \quad \text{on } \Omega \\
\dot{\mathbf{d}}(\mathbf{x},0) = \mathbf{v}_0(\mathbf{x}), & \quad \text{in } \Omega
\end{aligned} \quad (4)$$

where $T$ is the total time, $\zeta$ the space of permissible displacement, $\boldsymbol{\sigma}$ the stress tensor, $\mathbf{F}_b$ the body force field, and $\mathbf{n}$ the unit outward normal vector to $\Gamma_F$.

The direct integration method is used to solve the above dynamic problem, and the whole time domain, $\{t_i\}_{i=0}^{N_t}$, is divided into multiple time intervals, where $t_i$ is the $i^{\text{th}}$ time step, $N_t$ the number of time steps. Thus, the time-discrete form of the dynamic problem in Eq.(1) is obtained as follows:

$$\mathbf{M}\ddot{\mathbf{d}}_i + \mathbf{C}\dot{\mathbf{d}}_i + \mathbf{K}\mathbf{d}_i = \mathbf{f}_i, \quad i = 0,...,N_t \quad (5)$$

where $\mathbf{f}_i$ is the force vector at time step $i$ and $\mathbf{d}_i$, the displacement, $\dot{\mathbf{d}}_i$, the velocity, $\ddot{\mathbf{d}}_i$, the acceleration vector at $i^{\text{th}}$ time step, respectively.

The mass matrix is obtained by interpolating the density variables $b_e$ for the volume fraction. In order to obtain a clear topology optimization result, the threshold projection function is applied for volume interpolation [34], described as follows:

$$\mathbf{V}_e(b_e) = \frac{\tanh(\chi\eta) + \tanh(\chi(b_e - \eta))}{\tanh(\chi\eta) + \tanh(\chi(1-\eta))} \tag{6}$$

where $\chi$ is the control parameter of the projection and $\eta$ the threshold density.

The stiffness matrix is derived by interpolation of the density variables $b_e$ for the stiffness of each element, which is based on the RAMP function [35, 36].

$$\mathbf{E}_e(b_e) = \frac{V_e(b_e)}{1 + \kappa(1 - V_e(b_e))} \tag{7}$$

where $\kappa$ is the penalty parameter.

Given that the numerical instability occurs when $b_e \to 0$, the interpolation functions of the volume and the stiffness are defined as follows:

$$\overline{\mathbf{V}}_e(b_e) = \varepsilon + (1-\varepsilon)\mathbf{V}_e(b_e)$$
$$\overline{\mathbf{E}}_e(b_e) = \varepsilon + (1-\varepsilon)\mathbf{E}_e(b_e) \tag{8}$$

where $\varepsilon$ is an Ersatz parameter, $\varepsilon \ll 1$. Finally, the corresponding damping matrix is obtained according to Eq.(2).

According to the above discretization of displacement fields, the discretized dynamic topology optimization problem is mathematically represented as follows:

$$\min_{b \in [0,1]^N} f(b, \mathbf{d}_0, ..., \mathbf{d}_{N_t})$$
$$\text{s.t. } g(b) = v^T \overline{\mathbf{V}}_e(b_e) - V \leq 0 \tag{9}$$
$$\text{with: } \mathbf{M}\ddot{\mathbf{d}}_i + \mathbf{C}\dot{\mathbf{d}}_i + \mathbf{K}\mathbf{d}_i = \mathbf{f}_i, \ i = 0, ..., N_t$$

where $v$ is the volume vector of the finite elements, $V$ the volume fraction of the structure defined by the designer.

The general form of the objective function, which can be used with any optimization objective, is shown in Eq.(9). In this study, the corresponding objective function can be written as, to minimize the mean dynamic compliance,

$$f(b, \mathbf{d}_0, ..., \mathbf{d}_{N_t}) = \frac{1}{N_t} \sum_{t=0}^{N_t} \mathbf{f}_i^T \mathbf{d}_i \tag{10}$$

for mean strain energy,

$$f(b, \mathbf{d}_0, ..., \mathbf{d}_{N_t}) = \frac{1}{2N_t} \sum_{t=0}^{N_t} \mathbf{d}_i^T \mathbf{K} \mathbf{d}_i \tag{11}$$

to minimize the square of the displacement at a target degree of freedom,

$$f(b, \mathbf{d}_0, ..., \mathbf{d}_{N_t}) = \frac{1}{N_t} \sum_{t=0}^{N_t} (\mathbf{L}^T \mathbf{d}_i)^2. \tag{12}$$

### *2.2. Full-order model for adjoint sensitivity analysis*

In this study, sensitivity analysis is indispensable for the topology optimization presented in the previous section, and the adjoint method is adopted for sensitivity analysis to avoid direct derivation of state variables in each optimization design. In addition, the discrete approach is used, where state and time variables are previously discretized, without the consistency errors compared to the continuum method. It should be noted that the derivation of the discrete adjoint sensitivity in this study is closely related to the given integration method and the specific model-order reduction method. It is therefore necessary to explain the direct integration method used, called HHT-$\alpha$ [37].

The HHT-$\alpha$ method, generalizing the Newmark-$\beta$ method [38], introduces a new parameter $\alpha$ to modify the governing equation in Eq. (5) as follows:

$$\mathbf{M}\ddot{\mathbf{d}}_i + (1-\alpha)\mathbf{C}\dot{\mathbf{d}}_i + \alpha \mathbf{C}\dot{\mathbf{d}}_{i-1} + (1-\alpha)\mathbf{K}\mathbf{d}_i + \alpha \mathbf{K}\mathbf{d}_{i-1} = (1-\alpha)\mathbf{f}_i + \alpha \mathbf{f}_{i-1},$$

$$i = 1, ..., N_t. \tag{13}$$

Then, the Newmark-$\beta$ finite difference relations are assumed in the time domain of $(i-1) \sim i$:

$$\begin{aligned}\dot{\mathbf{d}}_i &= \dot{\mathbf{d}}_{i-1} + [(1-\delta)\ddot{\mathbf{d}}_{i-1} + \delta\ddot{\mathbf{d}}_i]\Delta t, \\ \mathbf{d}_i &= \mathbf{d}_{i-1} + \dot{\mathbf{d}}_{i-1}\Delta t + [(\tfrac{1}{2}-\beta)\ddot{\mathbf{d}}_{i-1} + \beta\ddot{\mathbf{d}}_i]\Delta t^2\end{aligned} \quad (14)$$

where $\beta$ and $\delta$ are parameters determined according to the requirements of integration accuracy and stability. The parameters $\alpha$, $\beta$, and $\delta$ are set in Eq.(15) to guarantee that the HHT-$\alpha$ method is at least second-order accurate and unconditionally stable [37]:

$$\begin{aligned}0 &\le \alpha \le \frac{1}{3}, \\ \beta &= \frac{(1+\alpha)^2}{4}, \\ \delta &= \frac{(1+2\alpha)}{2}.\end{aligned} \quad (15)$$

Substituting Eq.(14) into Eq.(5) gives the two-stage recursive formula in residual form:

$$\mathbf{R}_i = \hat{\mathbf{M}}_1\ddot{\mathbf{d}}_i + \hat{\mathbf{M}}_0\ddot{\mathbf{d}}_{i-1} + \hat{\mathbf{C}}_0\dot{\mathbf{d}}_{i-1} + \mathbf{K}\mathbf{d}_{i-1} - (1-\alpha)\mathbf{f}_i - \alpha\mathbf{f}_i = 0 \quad (16)$$

where

$$\begin{aligned}\hat{\mathbf{M}}_1 &= \mathbf{M} + (1-\alpha)\delta\mathbf{C}\Delta t + (1-\alpha)\beta\mathbf{K}\Delta t^2 \\ \hat{\mathbf{M}}_0 &= (1-\alpha)(1-\delta)\mathbf{C}\Delta t + (1-\alpha)(\tfrac{1}{2}-\beta)\mathbf{K}\Delta t^2 \\ \hat{\mathbf{C}}_0 &= \mathbf{C} + (1-\alpha)\mathbf{K}\Delta t\end{aligned} \quad (17)$$

The $\mathbf{d}_0$ and $\dot{\mathbf{d}}_0$ are given as initial conditions. First, calculate $\ddot{\mathbf{d}}_0$ as $\ddot{\mathbf{d}}_0 = \mathbf{M}^{-1}(\mathbf{f}_0 - \mathbf{C}\dot{\mathbf{d}}_0 - \mathbf{K}\mathbf{d}_0)$. Then calculate the nodal acceleration $\ddot{\mathbf{d}}_i$ in Eq. (14) and update $\mathbf{d}_i$ and $\dot{\mathbf{d}}_i$ according to Eq. (12) for each time step $i = 1,...,N_t$.

Following the HHT method described above, the derivation process of the sensitivity analysis is presented next. The sensitivity of the objective function $f(b, \mathbf{d}_0,...,\mathbf{d}_{N_t})$ to a design variable, $b_e$, is generally expressed as:

$$\frac{df}{db_e} = \frac{\partial f}{\partial b_e} + \sum_{i=0}^{N_t} \frac{\partial f}{\partial \mathbf{d}_i} \cdot \frac{\partial \mathbf{d}_i}{\partial b_e} \tag{18}$$

The adjoint method is used to avoid the direct derivation $\frac{\partial \mathbf{d}_i}{\partial b_e}$, which is computationally expensive. In this study, the recursive formula in residual form in Eq. (16) and the finite difference relationship based on the Newmark-$\beta$ method are taken as the adjoint term for the adjoint method [39]. In order to facilitate the deductions, Eq. (12) is rewritten in the form of residual as follows:

$$\begin{aligned}\mathbf{R}^{\mathbf{d}}_i &= -\dot{\mathbf{d}}_i + \dot{\mathbf{d}}_{i-1} + [(1-\delta)\ddot{\mathbf{d}}_{i-1} + \delta\ddot{\mathbf{d}}_i]\Delta t = 0 \\ \mathbf{R}^{\mathbf{v}}_i &= -\mathbf{d}_i + \mathbf{d}_{i-1} + \dot{\mathbf{d}}_{i-1}\Delta t + [(\tfrac{1}{2}-\beta)\ddot{\mathbf{d}}_{i-1} + \beta\ddot{\mathbf{d}}_i]\Delta t^2 = 0\end{aligned} \tag{19}$$

Then, the adjoint variables $\vartheta_i$, $\varsigma_i$, and $\tau_i$, $i = 0, \ldots, N_t$, are introduced and Eq. (18) is rewritten as

$$\begin{aligned}\frac{df}{db_e} &= \frac{\partial f}{\partial b_e} + \sum_{i=0}^{N_t} \frac{\partial f}{\partial \mathbf{d}_i} \cdot \frac{\partial \mathbf{d}_i}{\partial b_e} \\ &+ \sum_{i=0}^{N_t} \vartheta_i^T \left[ \frac{\partial \mathbf{R}_i}{\partial b_e} + \sum_{j=0}^{N_t} \left( \frac{\partial \mathbf{R}_i}{\partial \mathbf{d}_j} \cdot \frac{\partial \mathbf{d}_j}{\partial b_e} + \frac{\partial \mathbf{R}_i}{\partial \dot{\mathbf{d}}_j} \cdot \frac{\partial \dot{\mathbf{d}}_j}{\partial b_e} + \frac{\partial \mathbf{R}_i}{\partial \ddot{\mathbf{d}}_j} \cdot \frac{\partial \ddot{\mathbf{d}}_j}{\partial b_e} \right) \right] \\ &+ \sum_{i=1}^{N_t} \varsigma_i^T \left[ \frac{\partial \mathbf{R}^{\mathbf{d}}_i}{\partial b_e} + \sum_{j=0}^{N_t} \left( \frac{\partial \mathbf{R}^{\mathbf{d}}_i}{\partial \mathbf{d}_j} \cdot \frac{\partial \mathbf{d}_j}{\partial b_e} + \frac{\partial \mathbf{R}^{\mathbf{d}}_i}{\partial \dot{\mathbf{d}}_j} \cdot \frac{\partial \dot{\mathbf{d}}_j}{\partial b_e} + \frac{\partial \mathbf{R}^{\mathbf{d}}_i}{\partial \ddot{\mathbf{d}}_j} \cdot \frac{\partial \ddot{\mathbf{d}}_j}{\partial b_e} \right) \right] \\ &+ \sum_{i=1}^{N_t} \tau_i^T \left[ \frac{\partial \mathbf{R}^{\mathbf{v}}_i}{\partial b_e} + \sum_{j=0}^{N_t} \left( \frac{\partial \mathbf{R}^{\mathbf{v}}_i}{\partial \mathbf{d}_j} \cdot \frac{\partial \mathbf{d}_j}{\partial b_e} + \frac{\partial \mathbf{R}^{\mathbf{v}}_i}{\partial \dot{\mathbf{d}}_j} \cdot \frac{\partial \dot{\mathbf{d}}_j}{\partial b_e} + \frac{\partial \mathbf{R}^{\mathbf{v}}_i}{\partial \ddot{\mathbf{d}}_j} \cdot \frac{\partial \ddot{\mathbf{d}}_j}{\partial b_e} \right) \right]\end{aligned} \tag{20}$$

Assume that $\frac{\partial \mathbf{d}_0}{\partial z_e} = 0$ and $\frac{\partial \dot{\mathbf{d}}_0}{\partial z_e} = 0$, then simplify Eq. (18). Finally, the sensitivity of the objective can be obtained as:

$$\frac{df}{db_e} = \frac{\partial f}{\partial b_e} + \sum_{i=0}^{N_t} \boldsymbol{\vartheta}_i^T \frac{\partial \mathbf{R}_i}{\partial b_e} \tag{21}$$

Thus, the corresponding adjoint problem should be defined and the specific derivation process can be referred to [39]. The adjoint variable can be obtained by solving the following equations. This is the full-order model, which will then be reduced.

$$\boldsymbol{\varsigma}_{N_t} = \frac{\partial f}{\partial \mathbf{d}_{N_t}}, \ \boldsymbol{\tau}_{N_t} = 0, \ \hat{\mathbf{M}}_1 \boldsymbol{\vartheta}_{N_t} = -\beta \boldsymbol{\varsigma}_{N_t} \Delta t^2 - \delta \boldsymbol{\tau}_{N_t} \Delta t, \text{ for } i = N_t, \tag{22}$$

$$\begin{cases} \boldsymbol{\varsigma}_{i-1} = \dfrac{\partial f}{\partial \mathbf{d}_{i-1}} + \mathbf{K}\boldsymbol{\vartheta}_i + \boldsymbol{\varsigma}_i \\ \boldsymbol{\tau}_{i-1} = \hat{\mathbf{C}}_0 \boldsymbol{\vartheta}_i + \boldsymbol{\varsigma}_i \Delta t + \boldsymbol{\tau}_i \\ \hat{\mathbf{M}}_1 \boldsymbol{\vartheta}_{i-1} = \hat{\mathbf{M}}_0 \boldsymbol{\vartheta}_i - \left[\beta \boldsymbol{\varsigma}_{i-1} + \left(\dfrac{1}{2} - \beta\right)\boldsymbol{\varsigma}_i\right]\Delta t^2 - \left[\delta \boldsymbol{\tau}_{i-1} + (1-\delta)\boldsymbol{\tau}_i\right]\Delta t \end{cases}, \text{ for } i = 2,...,N_t - 1, \tag{23}$$

$$\mathbf{M}\boldsymbol{\vartheta}_0 = \hat{\mathbf{M}}_0 \boldsymbol{\vartheta}_1 - \left(\frac{1}{2} - \beta\right)\boldsymbol{\varsigma}_1 \Delta t^2 - (1-\delta)\boldsymbol{\tau}_1 \Delta t, \text{ for } i = 0. \tag{24}$$

Obviously, the dimension of the above adjoint problem scales with the DOF of the dynamic system, which can be extremely large when fine spatial discretization meshes are used. To deal with the complexity of the adjoint sensitivity analysis of the full-order model in Eqs. (20~22), the model order reduction technique is proposed next.

## 3. Fast adjoint sensitivity analysis based on reduced basis method

In this study, the reduced basis method aims to reduce the adjoint problem for sensitivity analysis, which requires repeated evaluations in topology optimization. The RBM consists of two phases, one being the offline phase where the main features of the adjoint problem described in Section 2.2 are captured based on the set of solutions for the full-order model, i.e. the truth solutions, for given density parameter values. Since the dimension $N_f$ of the adjoint problem solved to find the truth solution is large, the RBM approximates this full-order space with a $N_r$-dimensional low-order subspace (typically

$N_r \ll N_f$). In the online phase, the full-order adjoint problem is then mapped onto a $N_r$-dimensional subspace by Galerkin projection, allowing fast computation of the adjoint problem as $N_r \ll N_f$. That is, the online computational cost of the adjoint problem depends on $N_r$, rather than $N_f$. In general, an advantage is that the adjoint problem is solved faster in the online phase compared to the full-order model owing to the low-dimensional reduced space. Moreover, the basis vectors can express the reduced adjoint solution in full-order space because they are extracted from the full-order model. In many situations, fine mesh refinement is required to obtained accurate results, which increases the dimension of the full-order model, but has no evident effect on the online cost.

Numerically, the reduced adjoint solution in the time moment $\mathbf{S}^{N_t}$ in Eqs. (22)~(24) can be expressed as a linear combination of multiple basis functions, as follows:

$$\mathbf{S}^{N_t}(b) = \sum_{i=1}^{N_r} a_i^{N_t}(b) \mathbf{r}_i \qquad (25)$$

where $\mathbf{r} = \{\mathbf{r}_i, i=1,...,N_r\}$ are the reduced basis vectors, which are the dominant mode of truth solutions of the given parameters $b$, $a$ the dynamics of the generalized coordinates for the reduced system. The solution of the $\mathbf{r}$ and $a$ is presented in the following sections and the result is the reduced adjoint sensitivity analysis.

### 3.1. Offline stage: construction of the basis functions

While there are several strategies for generating reduced basis spaces, the POD-greedy algorithm is one of the most widely used [26]. POD lies in compressing the truth solutions at all sampling points of the specific parameter space, preserving only the dominant information, as explained in Algorithm 1. Further, the greedy algorithm is combined to generate the reduced basis space by iteratively adding a new basis vector and enriching the overall basis set.

**Algorithm 1**: POD algorithm, POD($\vartheta$, $n_{basis}$)

**Input**: Snapshots $\vartheta = \{\vartheta^1,...,\vartheta^{N_t}\} \in \mathbb{R}^{N_f \times N_t}$, number of the basis functions $n_{basis}$

**Output**: $\mathbf{r}_i \in \mathbb{R}^{N_f \times n_{basis}}$

**1** Make a Singular Value Decomposition for the snapshots, $\vartheta = \mathbf{MSV}$.

**2** Select the first $n_{basis}$ functions of $\mathbf{M}$, $\mathbf{r}_i = \mathbf{M}(:,1:n_{basis})$.

In general, the error estimator should be utilized for the POD-greedy algorithm, which avoids computing the full-order model in the overall sample space. Thus, the efficiency of basis function generation in POD-greedy can be improved in the offline phase. In this study, a strategy for generating reduced basis functions in POD-greedy is developed using the ML-based estimator, where the details of the proposed error estimator are detailed in Section 3.3. The generation of reduced basis vectors in the proposed approach is explained in Fig. 1.

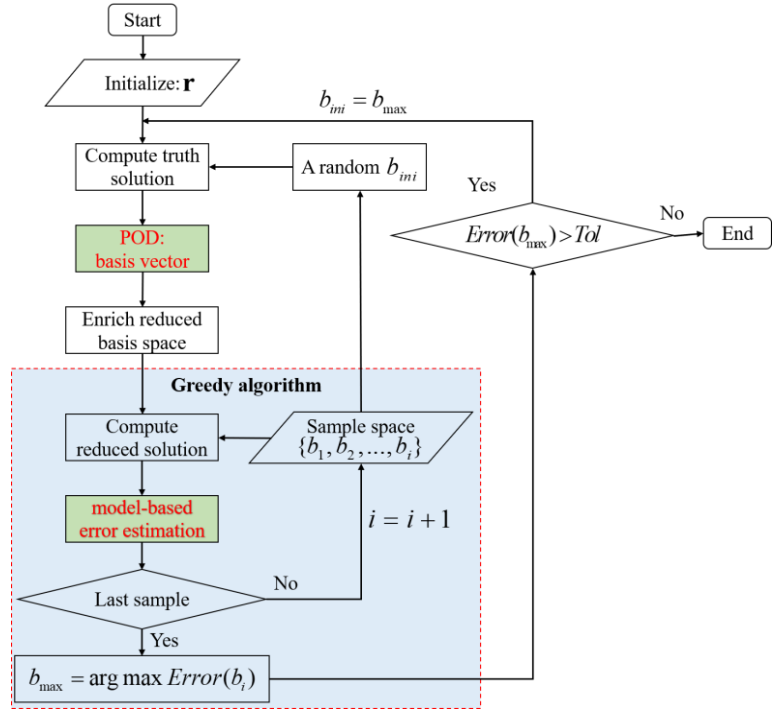

Fig. 1 The schematic representation of the proposed POD-greedy algorithm

In the offline phase, the reduced basis functions are generated as follows:

A parameter set $b_i$ is randomly selected from the sample space, and the true solution of the adjoint problem is computed according to Eq. (23). And the snapshot is obtained as $\vartheta(b_i) = \{\vartheta^n(b_i)\}_{n=0}^{N_t}$, i.e., $Snapshot_{b_i}$.

Calculate the reduced basis vector $r_i$ as in Algorithm 1, i.e., $r_i = \text{POD}(Snopshot_{b_i}, 1)$, and enrich the reduced basis space, $\mathbf{r} = [\mathbf{r}, r_i]$, where $r_i$ is the first POD mode of $\{\vartheta^n(b_i)\}_{n=0}^{N_t}$.

For all members of $b_i$, solve Eq. (25) to obtain the approximate solution $\mathbf{S}(b_i) = \{\mathbf{S}^n(b_i)\}_{n=0}^{N_t}$. The computational details of the reduced adjoint problem will be presented in the Section 3.2.

Perform the greedy algorithm for all $b_i \in D_h$ and find $b_{\max} = \arg\max_{b \in D_h} Error(b_i)$, calculated by the proposed model-based error estimator, which will be explained in Section 3.3.

Finally, it is judged whether the $Error(b_{\max})$ exceeds the prescribed threshold $Tol$. If the convergence condition is satisfied, the updated reduced basis space will be generated according to the updated $b_{\max}$, otherwise the cycle is exited and the final reduced basis space $\mathbf{r}$ is obtained.

### *3.2. Reduced-order model by Galerkin projection*

To facilitate the subsequent derivation of the reduced-order model and to simplify the notation, the parameter dependence on $b$ is removed, which can be contained due to the affine nature of the system [40]. The reduced basis space $\mathbf{r} \in \mathbb{R}^{N_r \times N_f}$ is determined in the offline phase described in Section 3.1 and the Galerkin projection is performed next. By substituting the ansatz formula (25) into the full-order model in Eqs. (22)~(24), the reduced-order model with the generalized coordinates $\mathbf{a}$ as state variables is obtained, i.e. the reduced-order adjoint problem, which is derived as follows:

$$\varsigma_{N_t} = \frac{\partial f}{\partial \mathbf{d}_{N_t}}, \tau_{N_t} = 0, \hat{\mathbf{M}}_1^r \mathbf{a}_{N_t} = \mathbf{r}^T(-\beta\varsigma_{N_t}\Delta t^2 - \delta\tau_{N_t}\Delta t), \text{ for } i = N_t, \tag{26}$$

$$\begin{cases} \varsigma_{i-1} = \dfrac{\partial f}{\partial \mathbf{d}_{i-1}} + \mathbf{K}^r \mathbf{a}_i + \varsigma_i \\ \tau_{i-1} = \hat{\mathbf{C}}_0^{\ r} \mathbf{a}_i + \varsigma_i \Delta t + \tau_i \\ \hat{\mathbf{M}}_1^r \mathbf{a}_{i-1} = \hat{\mathbf{M}}_0^{\ r} \mathbf{a}_i - \mathbf{r}^T \left[ \beta\varsigma_{i-1} + \left(\dfrac{1}{2}-\beta\right)\varsigma_i \right]\Delta t^2 - \left[\delta\tau_{i-1} + (1-\delta)\tau_i\right]\Delta t \end{cases}, \tag{27}$$

$$\text{for } i = 2, \ldots, N_t - 1,$$

$$\mathbf{M}^r \mathbf{a}_0 = \hat{\mathbf{M}}_0^{\ r}\mathbf{a}_1 - \mathbf{r}^T\left[\left(\frac{1}{2}-\beta\right)\varsigma_1\Delta t^2(1-\delta)\tau_1\Delta t\right], \text{ for } i = 0. \tag{28}$$

where $\hat{\mathbf{M}}_1^r = \mathbf{r}^T\hat{\mathbf{M}}_1\mathbf{r}$, $\hat{\mathbf{M}}_0^{\ r} = \mathbf{r}^T\hat{\mathbf{M}}_0\mathbf{r}$, $\mathbf{M}^r = \mathbf{r}^T\mathbf{M}\mathbf{r}$, $\mathbf{K}^r = \mathbf{K}\mathbf{r}$, $\hat{\mathbf{C}}_0^{\ r} = \hat{\mathbf{C}}_0\mathbf{r}$, $\mathbf{a} = [a_1, \ldots, a_{N_r}]^T$. Obviously, $\hat{\mathbf{M}}_1^r, \hat{\mathbf{M}}_0^{\ r}, \mathbf{M}^r \in \mathbb{R}^{N_r \times N_r}$ and $\mathbf{K}^r, \hat{\mathbf{C}}_0^{\ r} \in \mathbb{R}^{N_r \times N_t}$ are both of low dimension. The computational efficiency of the above reduced-order problem is much higher than that of the full-order model.

Then, the generalized coordinate $\mathbf{a}$ and reduced basis vectors constructed in Section 3.1 are substituted into Eq. (25) to obtain the approximate solution of the studied adjoint problem, and the sensitivity analysis of topology optimization is efficiently performed.

Finally, in order to improve the stability of the method, an adaptive adjoint problem computation strategy is introduced in the dynamic topology optimization process (online stage). Before computing each adjoint problem, if

$$countif(|b_i^{(k)} - b_i^{(k-1)}| \geq \varepsilon_1) \leq n \times \varepsilon_2 : i = 1, \ldots, n \tag{29}$$

the given reduced order model is computed, otherwise the full-order adjoint sensitivity analysis is executed. The convergence parameter $\varepsilon_1$ is an arbitrary value between 0 and 1 since $b_i$ has a range of [0,1]. $\varepsilon_2$ is a certain percentage.

*3.3. Model-based error estimation*

As described above, the reduced basis solution, $\mathbf{S}^{N_t}$, is calculated by the reduced basis method from Eqs. (26)~(28). Obviously, $\mathbf{S}^{N_t}$ satisfies the reduced-order model, yet not quite the full-order model in Eqs. (22)~(24). The full-order dynamics are formulated by substituting $\mathbf{S}^{N_t}$ into the full-order model, in which a residual $\mathbf{R}^{N_t} \in \mathbb{R}^{N_f}$ is present at each time step, described as follows:

$$\varsigma_{N_t} = \frac{\partial f}{\partial \mathbf{d}_{N_t}}, \tau_{N_t} = 0, \hat{\mathbf{M}}_1 \mathbf{S}_{N_t} = -\beta \varsigma_{N_t} \Delta t^2 - \delta \tau_{N_t} \Delta t + \mathbf{R}_{N_t}, \text{ for } i = N_t, \tag{30}$$

$$\begin{cases} \varsigma_{i-1} = \dfrac{\partial f}{\partial \mathbf{d}_{i-1}} + \mathbf{K}\mathbf{S}_i + \varsigma_i \\ \tau_{i-1} = \hat{\mathbf{C}}_0 \mathbf{S}_i + \varsigma_i \Delta t + \tau_i \\ \hat{\mathbf{M}}_1 \mathbf{S}_{i-1} = \hat{\mathbf{M}}_0 \mathbf{S}_i - \left[\beta \varsigma_{i-1} + \left(\dfrac{1}{2} - \beta\right)\varsigma_i\right] \Delta t^2 - \left[\delta \tau_{i-1} + (1-\delta)\tau_i\right]\Delta t + \mathbf{R}_i \end{cases}, \tag{31}$$

$$\text{for } i = 2, \ldots, N_t - 1,$$

$$\mathbf{M}\mathbf{S}_0 = \hat{\mathbf{M}}_0 \mathbf{S}_1 - \left(\frac{1}{2} - \beta\right)\varsigma_1 \Delta t^2 - (1-\delta)\tau_1 \Delta t + \mathbf{R}_0, \text{ for } i = 0. \tag{32}$$

Define the reduced-order error by $\mathbf{e}^{N_t} = \vartheta^{N_t} - \mathbf{S}^{N_t}$, and substitute it into the full-order dynamics above. This gives the error dynamics:

$$\varsigma_{N_t} = \frac{\partial f}{\partial \mathbf{d}_{N_t}}, \tau_{N_t} = 0, \hat{\mathbf{M}}_1 \mathbf{e}_{N_t} = -\beta \varsigma_{N_t} \Delta t^2 - \delta \tau_{N_t} \Delta t - \mathbf{R}_{N_t}, \text{ for } i = N_t, \tag{33}$$

$$\begin{cases} \varsigma_{i-1} = \dfrac{\partial f}{\partial \mathbf{d}_{i-1}} + \mathbf{K}\mathbf{e}_i + \varsigma_i \\ \tau_{i-1} = \hat{\mathbf{C}}_0 \mathbf{e}_i + \varsigma_i \Delta t + \tau_i \\ \hat{\mathbf{M}}_1 \mathbf{e}_{i-1} = \hat{\mathbf{M}}_0 \mathbf{e}_i - \left[\beta \varsigma_{i-1} + \left(\dfrac{1}{2} - \beta\right)\varsigma_i\right] \Delta t^2 - \left[\delta \tau_{i-1} + (1-\delta)\tau_i\right]\Delta t - \mathbf{R}_i \end{cases}, \tag{34}$$

$$\text{for } i = 2, \ldots, N_t - 1, ,$$

$$\mathbf{Me}_0 = \hat{\mathbf{M}}_0 \mathbf{e}_1 - \left(\frac{1}{2} - \beta\right)\varsigma_1 \Delta t^2 - (1-\delta)\tau_1 \Delta t - \mathbf{R}_1, \text{ for } i = 0. \tag{35}$$

Since the above error dynamics have the dimension of the full-order model, it is clear that they cannot be used to efficiently determine the estimation error. Thus, additional requirements for error estimation are imposed by the need for efficient and accurate error bounds.

In accordance with the $\ell_2$-gain theory, the reduction of the system error should be represented as:

$$\|\mathbf{e}\|_2 \leq \lambda \|\mathbf{R}\|_2 \tag{36}$$

where $\lambda$ is the $\ell_2$-gain of the system from input $b$ to output $\vartheta$ and $\|\mathbf{R}\|_2$ is the two-norm of the residual over the overall DOF at each instant.

To order to efficiently solve the norm involved in the error estimation in Eq. (36), Abbasi et al. [26] propose a concise approach based on Euclidean distance. For this method, in the offline phase, the required $\lambda$ corresponding to each parameter in the discrete parameter domain $D_{offline} \subset D_b$, is calculated and stored. In the online phase, the error estimation is performed and the $\lambda$ is determined according to the Euclidean distance between the new online parameter $b^e$, which does not lie in $D_{offline}$, and each offline parameter instead of solving directly. The closest parameter $\bar{b}^e$ in the $D_{offline}$ is found, which is defined in Eq. (43). Since only the norms corresponding to the parameters of $D_{offline}$ are stored, the efficiency of the error estimation is improved and the memory requirement is not increased.

$$\bar{b}^e = \arg \min_{b \in D_{offline}} \|b - b^e\| \tag{37}$$

The essence of this method is to establish the relationship between error and residual in advance based on the offline discrete parameter domain, and the error corresponding to the new parameter is obtained directly according to the previously stored relationship. The disadvantage of this method is that the requirements for the discrete parameter domain are relatively high. Since new parameters need to find the corresponding similar parameters in the offline parameter domain, intensive offline parameter partitioning must be ensured to achieve more accurate error prediction. Moreover, the new parameters and offline parameters only use the Euclidean distance to judge the approximation between parameters, which will lead to multiple parameters corresponding to the same distance, causing instability in error estimation. In topology optimization, each cell corresponds to a density parameter, which makes the problem studied have a large number of design parameters. The increasing number of parameters aggravates the shortcomings of this error estimation method. On the one hand, it requires a large number of offline samples, which makes it prohibitively expensive, especially for dynamic problems. On the other hand, it is more difficult to search for the nearest parameter in the offline parameter domain for new parameters.

Given the above bottlenecks, standard NN model is developed in this study to establish the relationship between the true error and residual. As we all know, NN model can deal with various complex modelling problems and has the ability to achieve accurate prediction based on a small number of samples. It is applied here to establish a model-based error estimator from modelling to prediction. The methodology and various components of the error estimator are described in the remainder of this section.

Among various types of NN models, the FNN has been widely used due to its excellent modelling capabilities and ease of implementation. In this study, a simple multi-layer feed-forward network including input, hidden and output layers, as shown in Fig. 2, is developed and trained by error back-propagation. By minimizing the square of the error of the network, the weights and thresholds of the neural network are continuously updated by the steepest descent method based on learning roles through error back-propagation.

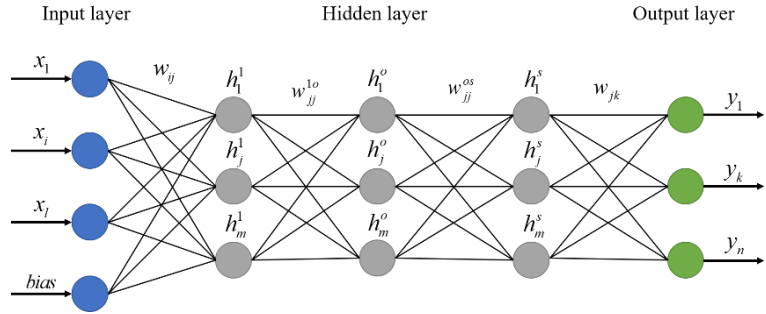

Fig. 2 The construction of the present FNN

## 4. Numerical case studies

In this section, numerical examples are discussed to substantiate the feasibility and effectiveness of the proposed method. The section consists of three subsections and is organized as follows: In Section 4.1, two 2D benchmark examples subjected to different periodic loads are demonstrated in terms of the minimum mean dynamic compliance and minimum squared displacement at a target degree of freedom problems. In Section 4.2, the building structure is studied with regard to ground excitation, and the acceleration effect of the reduced-order model is investigated to show its influences on the 2D structure with variant scales. In Section 4.3, the 3D numerical studies are presented to further express the effectiveness and efficiency of the algorithm, including verifications on random loading condition and different structure scales. Also, $\alpha = 0.05$, $\beta = (1+\alpha)^2/4$, $\delta = (1+2\alpha)/2$ are given to guarantee at least second-order accurate and unconditional stable in the HHT-$\alpha$ method. In addition, the validity of the proposed error estimation method is verified in terms of accuracy and efficiency, which is illustrated by the following examples.

## 4.1. Beam structure design under periodic load

In this section, two 2D benchmark examples with a regular design domain are performed for two types of optimization problems - minimum squared displacement at a target degree of freedom and minimum mean dynamic compliance. The definition and illustration of all benchmark examples are outlined in the Fig. 3, including boundary conditions and external loading. These structures have size, $L = 4$ m, thickness, $h = 0.01$ m, and their material is steel with Young's modulus, $E_0 = 200$ GPa, Poisson's ratio, $\nu = 0.3$, and mass density, $\rho_0 = 7800$ kg/m³. To solve this dynamic problem, the number of time steps is set to 200 and simulation time to 0.05s.

The first example is a cantilever beam subjected to a half-cycle sinusoidal load applied at the centre of the free edge of the beam. For this example, the objective function is the mean square displacement at the DOF where the load is applied. The second example contains a square domain supported at the bottom and forced at the middle of the top with a constant load rotating at a given angular frequency, $\omega = 20\pi$ rad/s.

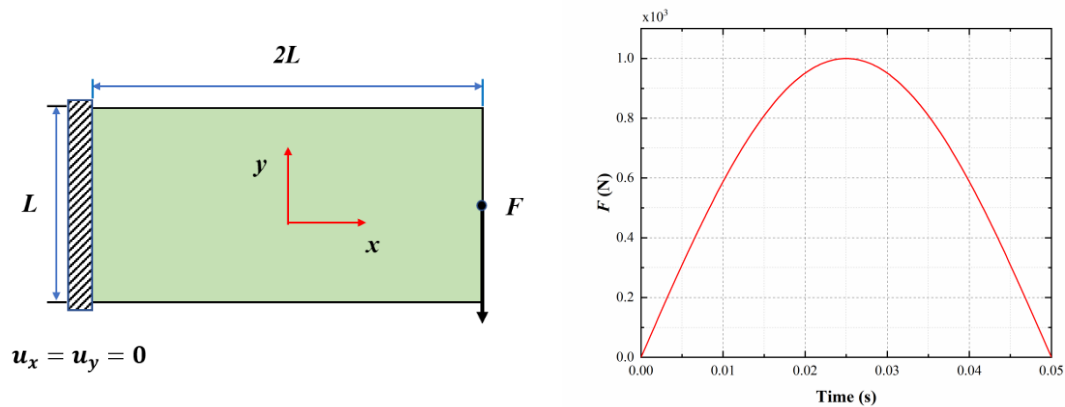

(a) Cantilever beam

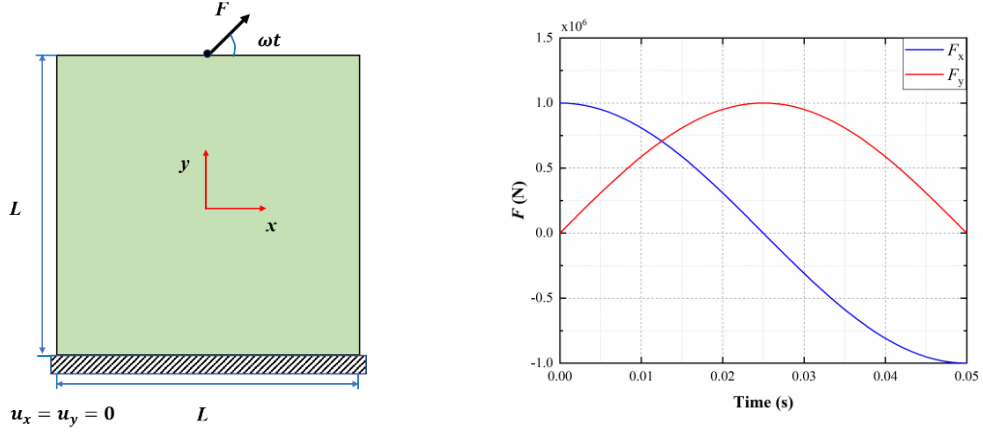

(b) Support structure
Fig. 3 The benchmark examples: boundary conditions and external load.

*5.1.1. Offline phase: Reduced basis vectors construction.*

In the offline phase, the POD-Greedy method is employed to generate an increasing sequence of reduced spaces $\{\mathbf{r}_{N_r}\}_{N_r>0}$. In these two examples, the discrete parameter domain $D_h$ contains 200 members, all from the full order models involved in different stages of the dynamic topology optimization process, distributed over the parameter domain $D$.

Throughout the greedy iteration process, a basis vector is incrementally added to update the reduced space at each iteration according to the error between the updated reduced-order model and the full-order model listed in $D_h$. The full-order adjoint solution that has the largest error with the reduced ones is used to enrich the reduced space. There, we define the $l_2$-norm of the error formula to evaluate the absolute error at each time interval of the dynamic adjoint sensitivity analysis. The error at each time interval was calculated by

$$error_{dy} = \left\| \mathbf{d}_{reduceed} - \mathbf{d}_{full} \right\|_2 \tag{38}$$

where $\mathbf{d}_{reduceed}$ is the approximate solution obtained by the proposed reduced-order method, $\mathbf{d}_{full}$ is the exact solution. After calculating the residual, we can build the expression in Eq. (36) for adjoint sensitivity analysis and take it as the relationship between residual and error. Then, the prediction problem can be written as

$$\{\mathbf{e}_1,...,\mathbf{e}_{N_t}\} = f(\{\mathbf{R}_1,...,\mathbf{R}_{N_t}\}) \tag{39}$$

The function in Eq. (39) is the error evaluator that is established by the FNN in Fig. 2. The trained error evaluator plays the role of snapshot selection in greedy programs. It is efficient to determine whether to add a new solution snapshot to the reduced space without time-consuming complete analysis. Fig. 4 shows the error evolution history of the incrementally updated reduced space. The abscissa represents the number of steps in the greedy iterative process, which is equal to the number of basis vectors, and the ordinate represents the error between the approximate adjoint solution and the exact adjoint solution as defined in Eq. (38). It can be seen that as the basis functions are enriched, the error gradually decreases and tends to be stable. Thus, the optimal set of basis functions is determined for the reduced-order model. For the cantilever beam design, the reduced order model is given with 11 basis functions and 37 for the support structure. Moreover, to demonstrate the accuracy of its error evaluation, we use root mean square error ($RMSE$) and coefficient of determination ($R^2$) as evaluation criteria, and the comparison between real error and evaluation error is also shown in the Fig. 4. It can be illustrated that the estimation error basically matches the true error throughout the entire time domain. The evaluation indicators ($RMSE$ and $R^2$) further validate that the present error estimation method can predict the true error well.

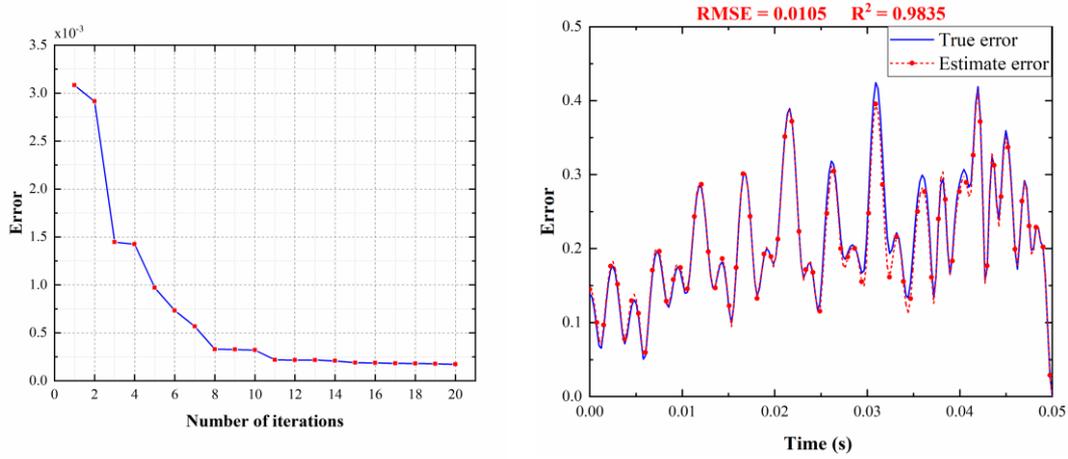

(a) Cantilever beam

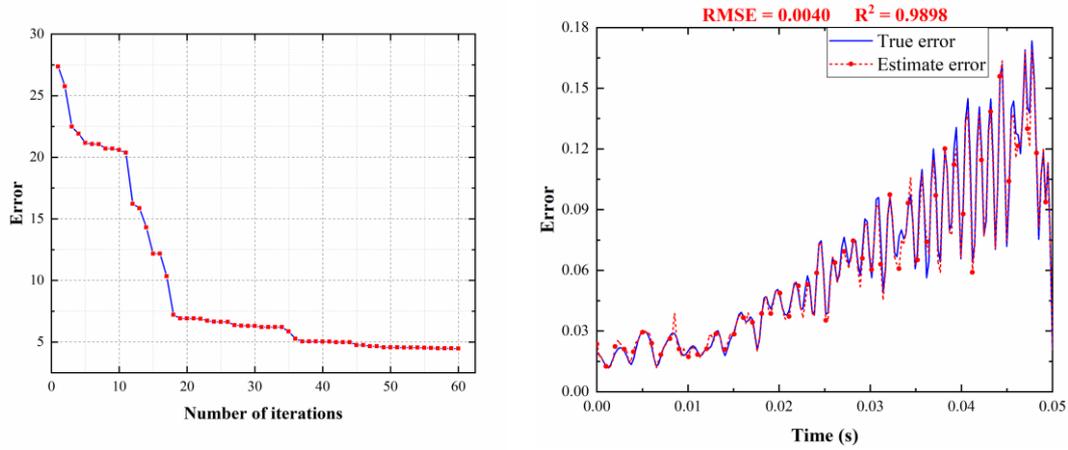

(b) Support structure

Fig. 4 Offline phase: Greedy iteration and error estimation.

*5.1.2. Online phase: Adaptive dynamic TO with reduced order model.*

The corresponding basis functions were obtained during the offline stage mentioned above, and then the topology optimization results could be quickly obtained by projection. The proposed method is mainly evaluated from two aspects: approximation accuracy and acceleration effect.

For the cantilever beam example, Fig. 5 shows the TO results to validate the accuracy, including the objective function, the topology configuration, and the comparison of the structural performance before and after optimization. The optimization objective of this problem is to minimize the vertical displacement of the loading point. Fig. 5 (a) and (b) show that the reduced-order model established by the present methodology achieves a topology configuration consistent with the full order model.

Besides, the evolution of the target displacement displayed in Fig. 5 (c) indicates that the maximum displacement at the interested point of the optimal structure is smaller than the initial structure.

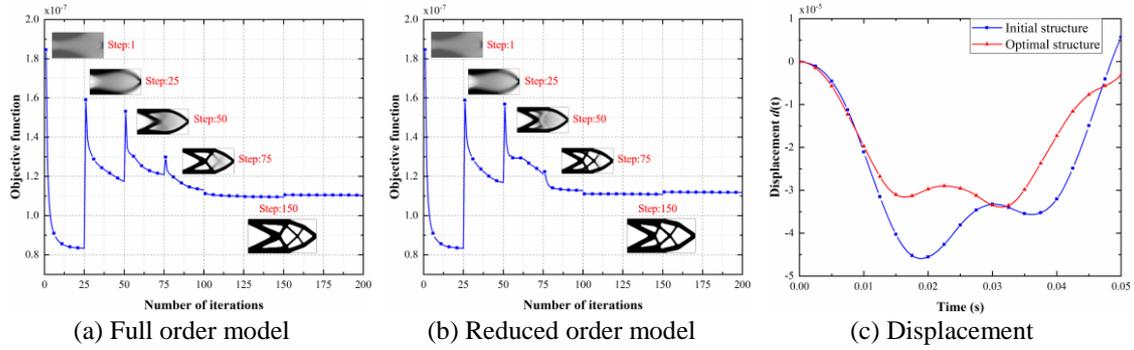

(a) Full order model　　　　(b) Reduced order model　　　　(c) Displacement
Fig. 5 Topology optimization results of cantilever beam.

For the support structure design, the objective function is to minimize the dynamic compliance of the structure, and the TO results are shown in Fig. 6. It can also be concluded that the proposed reduced-order model approximates well and achieves a topology configuration consistent with the full order model. In addition, the reduction in structure deformation compared to the initial structure indicates an improvement in the stiffness performance of the optimal structure.

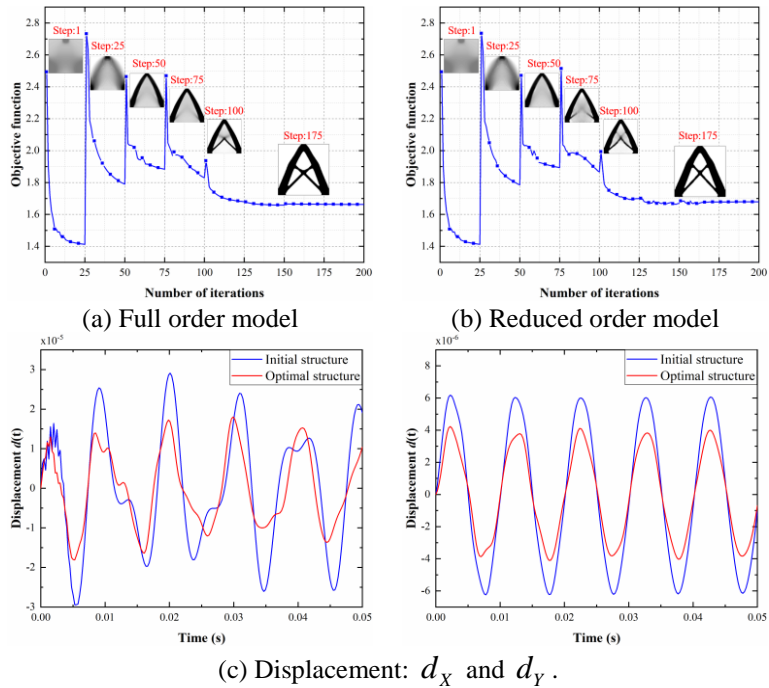

(a) Full order model　　　　(b) Reduced order model

(c) Displacement: $d_X$ and $d_Y$.
Fig. 6 Topology optimization results of support structure.

Next, the computational efficiency of the proposed reduced-order model is discussed. The cantilever beam domain is discretized into linear quadrilateral elements with 10368 elements and 10585 nodes, and the support structure finite element model with 10000 elements and 10201 nodes. Fig. 7 shows the computation time of the adjoint problem in each dynamic topology iteration. As an adaptive strategy is developed to dynamically select models for the computation of adjoint problems, the distribution of points in the figure clearly indicates its triggering time. Owing to significant structural changes in the early stage of optimization, the full-order model is employed repeatedly to ensure the accuracy of adjoint solution. In the later stage of the optimization, only the reduced-order model is used to satisfy the accuracy requirements well. Further, the computation time of the reduced-order model is significantly lower than that of the full-order model from a single iteration perspective. In summary, the efficiency of the proposed method has been verified.

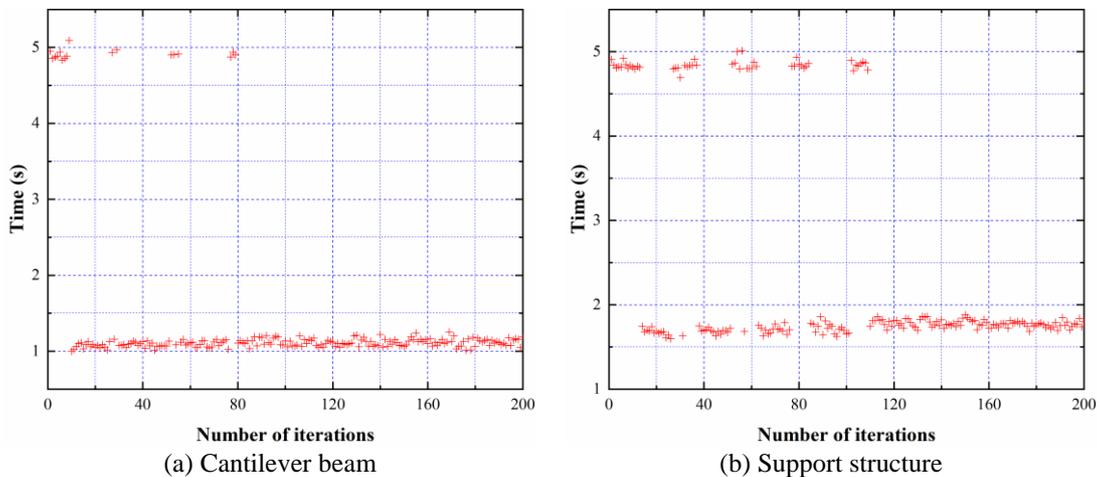

(a) Cantilever beam  (b) Support structure
Fig. 7 The calculation time of adjoint sensitivity analysis.

*4.2. Building design under ground excitation*

In this example we aim to minimize mean strain energy under ground excitation. The design domain and load profiles are shown in Fig. 8. This is a building design domain, the bottom end is fixed and the top is subject to the ground acceleration of a lumped mass point varying sinusoidally in time. The building has height $H = 75$ m, width $L = 30$ m, thickness $h = 1$ m, and a lumped mass of magnitude $M = 0.4 \times 10^6$ kg, as shown in Fig. 8 (a). The acceleration curve is shown in Fig. 8 (b), $a_g(t) = 5 \times \sin(2.5\pi \times t)$ m/s². In addition, the adjoint sensitivity analysis of the building is divided into 200 intervals within [0,4.8] sec. In this section, variant dimensions of the reduced-order models of the adjoint problem are performed to explore their effects on the optimal topology configuration, which further validates that the present reduced-order model has the characteristic of balancing efficiency and accuracy. Furthermore, the acceleration effect of the given reduced-order model on 2D problems with different scales was investigated.

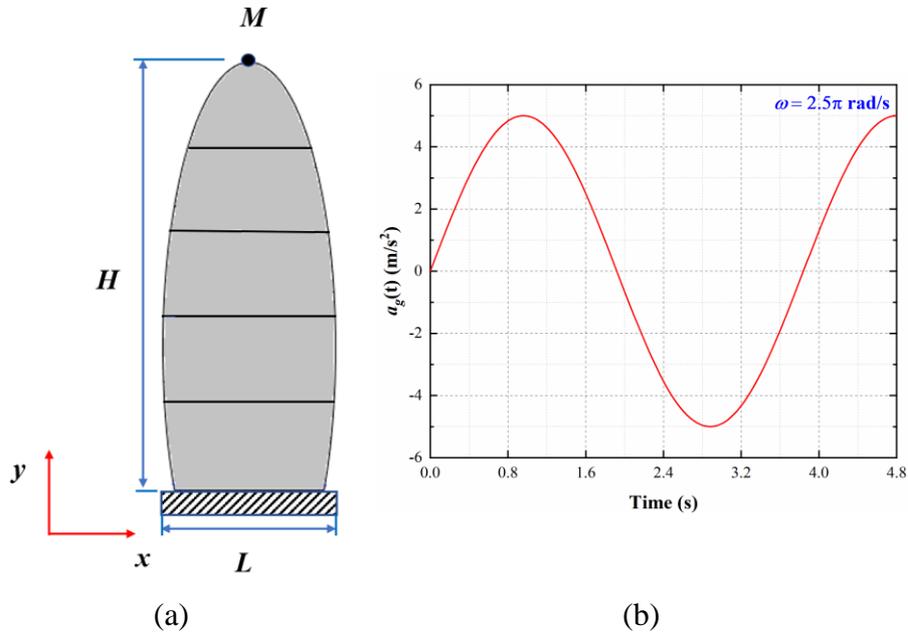

(a)                  (b)
Fig. 8 The building structure: (a) domain and boundary conditions and (b) applied acceleration.

*5.2.1. Dimension of reduced-order model*

The dimension of the reduced-order model directly affects the accuracy and efficiency of the sensitivity analysis. In this method, the number of basis functions is controlled by a given error threshold in the offline stage. Theoretically, if the set error is small enough, an optimal reduced order model in terms of accuracy can be obtained. However, excessive pursuit of a high accuracy approximation leads to a decrease in the computational efficiency of the reduced-order model. Therefore, several reduced-order models with different dimensions are discussed next.

Fig. 9 shows the construction process of the reduced basis functions in the offline phase. With the enrichment of the basis functions, the error between the reduced-order model and the full ones in the offline sample space gradually decreases. From the perspective of only pursuing the accuracy of the reduced order model, when the number of basis functions reaches 180, the error trend tends to stabilize and approaches 0. However, in order to simultaneously pursue high computational efficiency, this reduced order model is unreasonable. Furthermore, from the optimal topological configurations corresponding to different reduced-order models in Fig. 9, it can be seen that the optimal structure obtained by the reduced-order model composed of multiple basis functions has the main features consistent with the full-order model, see Step. 10; as the basis functions in the reduced-order model are continuously enriched, the ability to capture the features is enhanced, and the small features in the optimal structure are gradually captured; for this example, the corresponding reduced-order model with 40 basis functions can achieve the optimal structure consistent with the full-order model when the error between the models is reduced to a certain level.

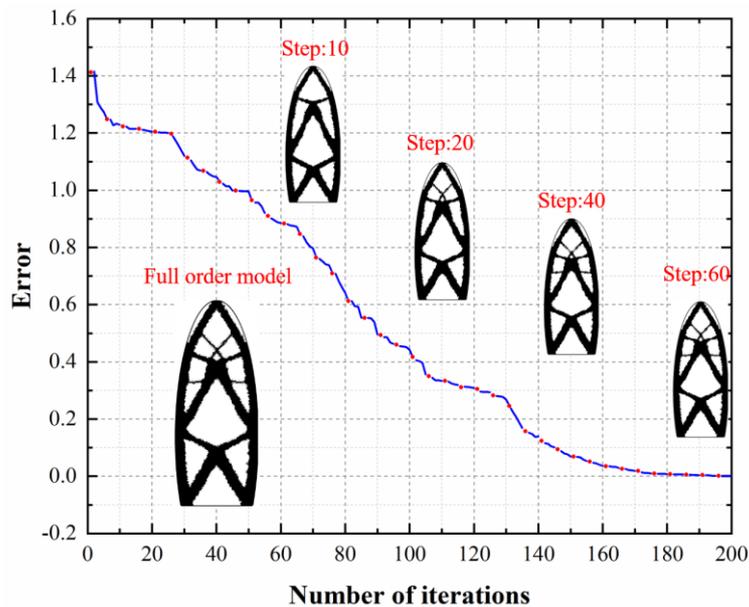

Fig. 9 Greedy iteration in Offline phase for building design.

The next step is to further test the representation of the computational efficiency of the optimization process for different reduced-order models. Fig. 10 (a) shows the execution of the reduced-order models of different dimensions after applying the adaptive strategy. The curve in the figure indicates the cumulative number of executions of the full-order model, with more triggers and a corresponding increase in the total iteration time. In the precondition, the same threshold parameters are used for each model. It is clear from the figure that, in general, the dimension of the reduced-order model is inversely proportional to the number of executions of the full-order model, because the higher the dimension, the higher the approximation accuracy of the model, and the higher the accuracy of the adjoint solutions in the optimization, the less likely it is that the execution conditions of the full-order model will be triggered. Fig. 10 (b) provides further statistics on the computational efficiency of the average per adjoint solution, expressed in terms of the speedup. And combining the results of the optimal topology of each model in Fig. 9, it can be concluded that although the speedup of the reduced-order model is better than that of the other models when taking 10 and 20 basis functions, its optimal structure lacks some small features and has insufficient accuracy compared to that of the full-order model. Furthermore, as the dimensionality of the reduced-order model increases, the approximation accuracy meets the requirements and an optimal structure consistent with the full-order model is obtained. However, the number of executions of the full-order model is not significantly reduced and the acceleration performance deteriorates. Thus, a reduced-order model with 40 basis functions is the optimal choice for this example.

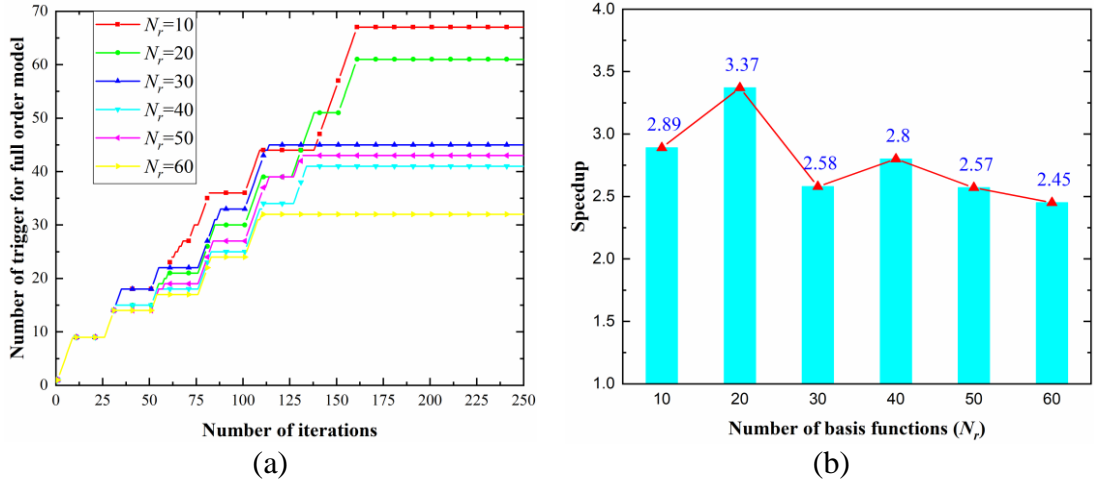

(a)                      (b)

Fig. 10 Adaptive model selection strategy:

(a) Execution of full-order model; (b) Speedup of adjoint problem for average per iteration.

*5.2.2. Scale of the adjoint problem*

Subsequently, building structures with variant scales are optimized with respect to the determined reduced space dimension, and the information of their finite element models is listed in Table 1.

Table 1 The finite element model information

| Building model | 1 | 2 | 3 | 4 | 5 | 6 |
|---|---|---|---|---|---|---|
| Node | 19707 | 39399 | 59007 | 79267 | 99438 | 198730 |
| Element | 10000 | 20000 | 30000 | 40000 | 50000 | 100000 |
| DOF | 39414 | 78798 | 118014 | 158534 | 198876 | 397460 |

In this subsection, we focus on exploring the application potential of the proposed reduced-order model to large-scale problems. The computational efficiency of using the reduced-order model and the full-order model is investigated comparatively in an adjoint sensitivity analysis. The speedup curves for applying the reduced-order model to solve adjoint problems of different scales are shown in Fig. 11. It can be seen that the speedup increases as the problem scale increases, and although the increasing trend is gradually flattened, the acceleration performance can still be maintained at a high level. Overall, this indicates that the proposed reduced-order model has the potential to be applied to large-scale problems.

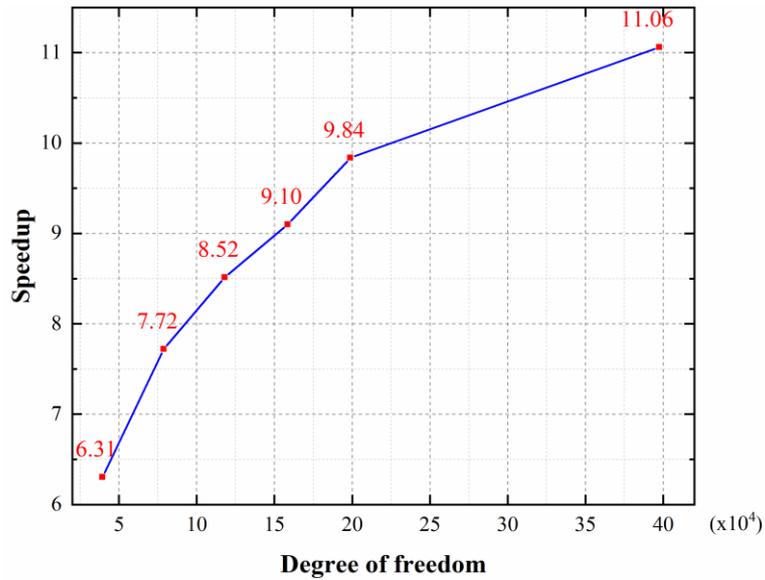

Fig. 11 Speedup of reduced-order model for adjoint sensitivity analysis.

Finally, dynamic topology optimization is performed based on a reduced adjoint sensitivity analysis approach, the results of which are shown in Fig. 12. The total strain energy of the structure is reduced after the optimization, and the displacement in the y-direction at the lumped mass point is discussed. As can be seen in Fig. 12 (b), the deformation of the structure is reduced and the dynamic stiffness performance is improved.

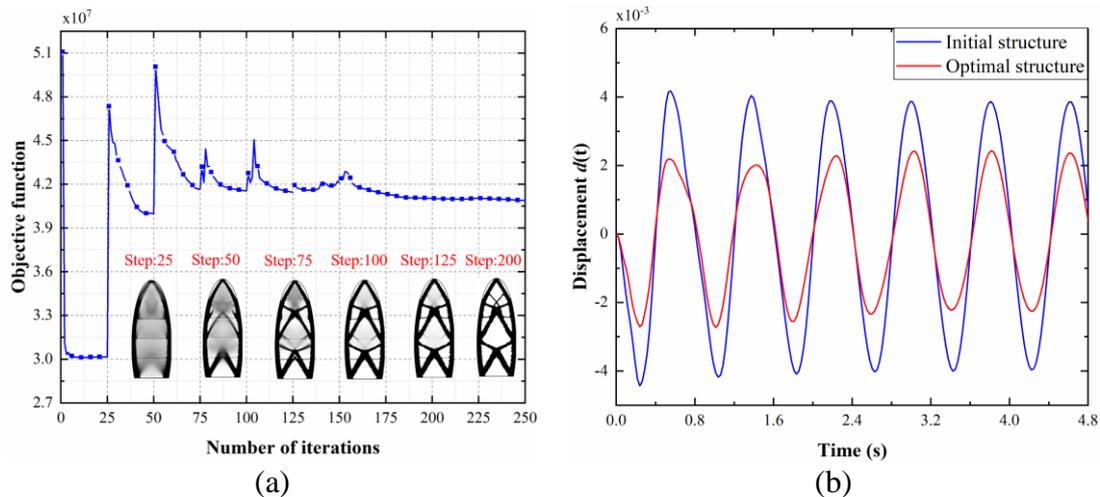

(a) (b)

Fig. 12 Topology optimization results: (a) Objective function; (b) Displacement.

### *4.3. 3D structure design under random excitation*

The 3D structure is further investigated to demonstrate the validity and feasibility of the developed reduced strategy. To this end, a wrench structure with height $H = 0.1$ m, width $W = 0.28$ m and thickness $Th = 0.02$ m is subjected to randomized dynamic force, as displayed in Fig. 13. The right hole wall surface of the wrench structure is constrained, and randomized dynamic loads are applied to the left half hole wall surface.

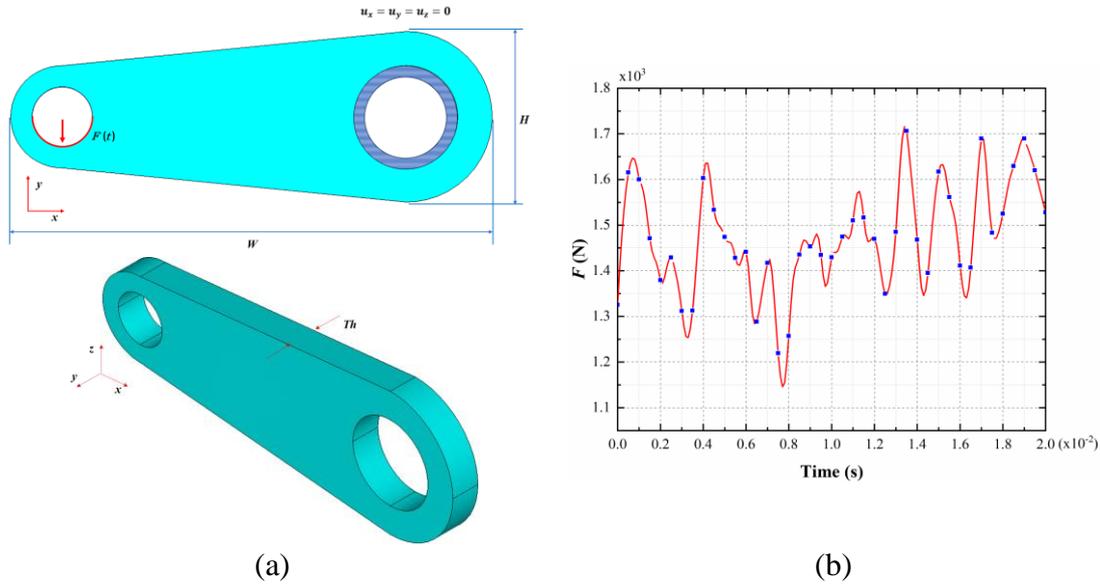

(a) (b)

Fig. 13 The wrench structure: (a) domain and boundary conditions and (b) applied load.

First, in the offline phase, the proposed model-assist greedy procedure is performed. Fig. 14 (a) shows the error evolution of the reduced-order model constructed in the greedy iteration, and it can be seen that the error gradually decreases as the basis functions are enriched, and the error stabilises at a very small value when the number reaches 10. Furthermore, the optimal structure obtained on the basis of the reduced-order model is checked for its accuracy. In this process, the error of the reduced-order model is obtained by the presented model-based estimator, which avoids solving the full-order adjoint problem and improves the efficiency. Fig. 14 (b) shows the results of the comparison between the true error and the estimated error. The evaluation indicators ($RMSE = 0.0013$, $R^2 = 0.9871$) validate that the present model can accurately predict the true error.

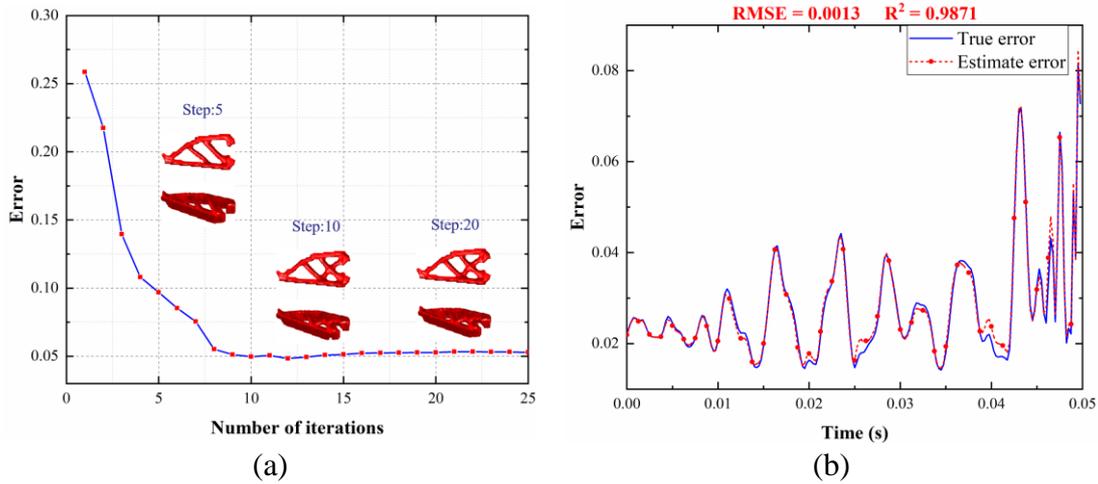

Fig. 14 Offline phase: (a) Greedy iteration; (b) Error estimation.

Then, in the online phase, topology optimization is carried out with the goal of minimizing dynamic compliance using adjoint methods for sensitivity analysis. At each iteration of solving the adjoint problem, a full-order model or a reduced-order model ($N_r = 10$) is selected based on the given threshold ($\varepsilon_1 = 0.1$, $\varepsilon_2 = 0.005$), which is a fast, easy-to-compute and reasonably accurate model built by Galerkin projection. Fig. 15 (a) shows the results of the adaptive strategy. It can be seen that the full-order model is used in the pre-optimization stage to ensure the optimization accuracy, and the reduced-order model is triggered in the post-optimization stage to significantly improve the computational efficiency while ensuring the accuracy. Besides, the accuracy is further verified by the topological configurations of some of the iterations in Table 2. In addition, the physical property of the optimal structure is verified in Fig. 15 (b), which displays the y-direction displacement of the load point. In conclusion, the significant reduction in structure deformation after optimization indicates an improvement in the stiffness property of the structure.

Table 2 Topology optimization results.

| Iteration | Full-order model | Reduced-order model |
|---|---|---|
| 1 | | |

| 5 | | | | |
| 10 | | | | |
| 20 | | | | |
| 50 | | | | |
| 100 | | | | |

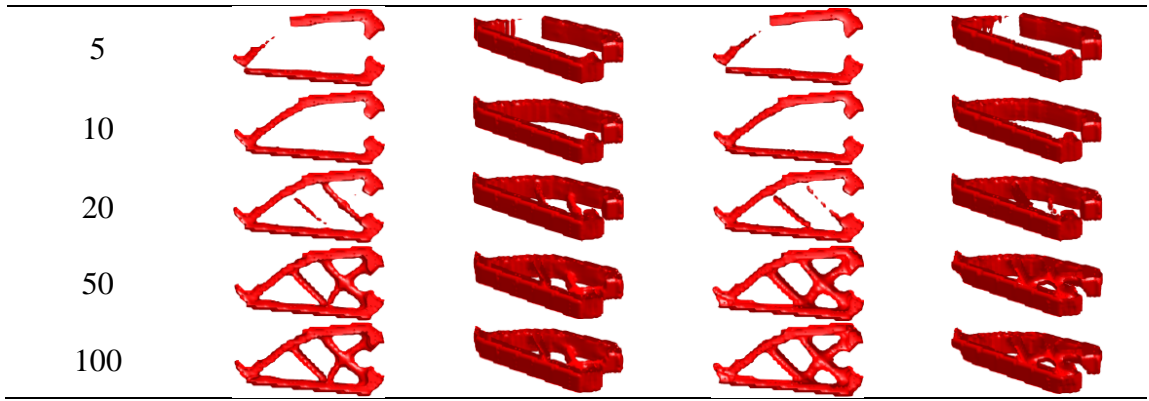

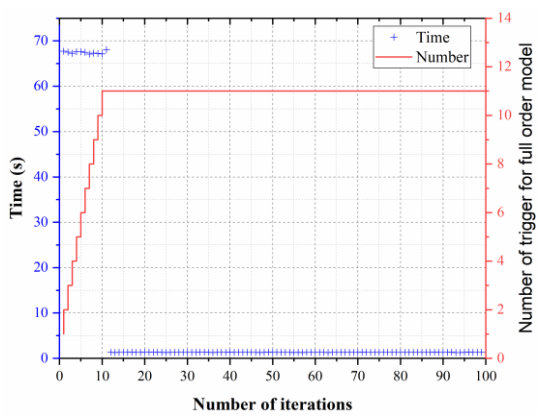

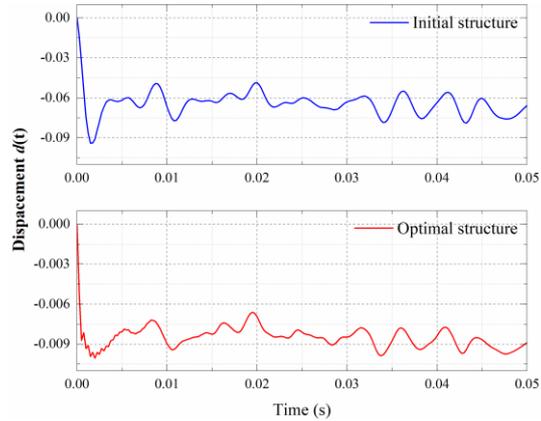

(a)            (b)

Fig. 15 The effect of reduced-order model:

(a) The calculation time of adjoint sensitivity analysis; (b) Displacement of the load point.

Finally, the acceleration performance of the reduced-order model for 3D structure design is discussed. And wrench structures of variant scales are optimized, which are listed in Table 3.

Table 3 The finite element model information

| Wrench model | 1 | 2 | 3 | 4 | 5 |
| --- | --- | --- | --- | --- | --- |
| Node | 3286 | 10408 | 19252 | 31930 | 45628 |
| Element | 540 | 1608 | 3528 | 6640 | 11100 |
| DOF | 9858 | 31224 | 57756 | 95790 | 136884 |

Fig. 16 shows the speedup curves for applying the reduced-order model to solve adjoint problems of variant scales. It can be concluded that the present methodology achieves better speedup when applied to 3D structures compared to 2D problems; as the problem scale increases, the speedup grows in a trend consistent with the 2D results, but ultimately maintains high acceleration performance.

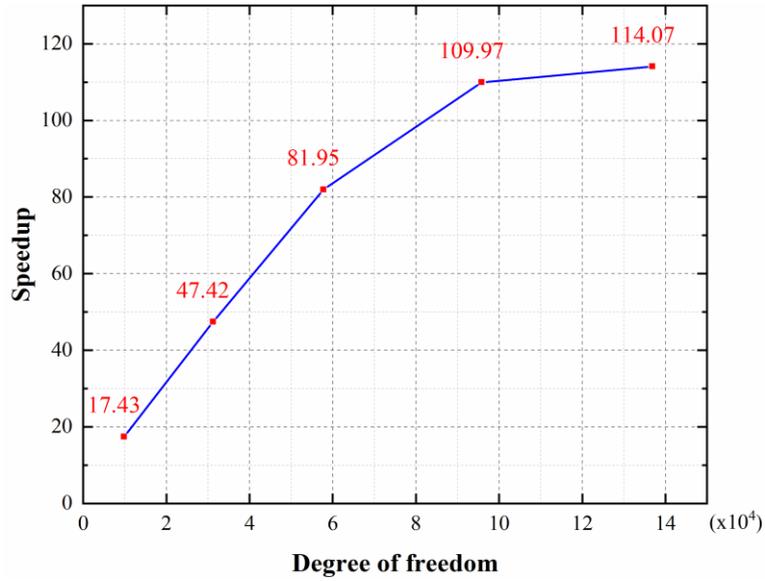

Fig. 16 Speedup of reduced-order model for adjoint sensitivity analysis

**Conclusions**

In this study, model-based error estimation is proposed and applied to the offline phase in RBM to speed up and facilitate the generation of reduced basis space. In particular, a reduced-order model method is developed for the adjoint sensitivity analysis of dynamic topology optimization with discrete sensitivity, which is rarely researched. In the proposed method, the main contributions can be summarized as follows:

(i) The RBM is introduced to establish the dynamic reduced-order model for the adjoint sensitivity analysis with discrete sensitivity;

(ii) Model-based error estimation is developed for the reduced-order model of the adjoint problem, and the typical FNN is constructed to establish the relationship between the residual and the true error;

(iii) The RBM is further developed with the model-based error estimation, which provides appropriate solutions of the adjoint problem combined with the greedy algorithm in the offline phase, thus speeding up and facilitating the generation of the reduced basis space.

Three examples with different initial design domains, load boundary conditions and model sizes are studied. It can be concluded that the proposed reduced-order model has the powerful ability to efficiently solve adjoint problems for discrete sensitivity analysis in dynamic topology optimization, while achieving excellent accuracy. Furthermore, the model-based error estimation brilliantly models the relationship between the residual error and the true error, which facilitates the RBM for the adjoint problem. Theoretically, the proposed reduced-order model method could be used in different integral solution frameworks of motion equations. Moreover, according to the complexity of different problems, the model-based error estimation should be pursued with different neural networks.

**Declaration of competing interest**

The authors declare that they have no known competing financial interests or personal relationships that could have appeared to influence the work reported in this paper.

**Acknowledgements**

We acknowledge the support provided by the Project of the National Natural Science Foundation of China (11702090) and Peacock Program for Overseas High-Level Talents Introduction of Shenzhen City (KQTD20200820113110016).